\documentclass{article}
\usepackage{amsmath,amssymb,amsthm}

\newcommand{\openR}{\mathbb{R}}

\newcommand{\pp}{\mathbb{P}}

\newcommand{\cov}{{\rm \mathbb{C}ov}}
\newcommand{\half}{\frac{1}{2}}
\newcommand{\ee}{\mathbb{E}\,}

\newcommand{\R}{\mathbb{R}\,}

\newcommand{\eps}{\varepsilon}

\newcommand{\bp}{\mbox{{\boldmath $\mathcal{P}$}}}
\newcommand{\bs}{\mbox{{\boldmath $\Sigma$}}}

\renewenvironment{proof}[1][\proofname]{\par \normalfont \trivlist
 \item[\hskip\labelsep\itshape #1]\ignorespaces
}{%
 \hspace*{\fill}$\Box$ \endtrivlist
} \renewcommand{\proofname}{{\bf Proof}}

\newtheorem{theorem}{Theorem}[section]
\newtheorem{proposition}[theorem]{Proposition}
\newtheorem{lemma}[theorem]{Lemma}
\newtheorem{corollary}[theorem]{Corollary}
\theoremstyle{definition}

\newtheorem{remark}[theorem]{Remark}
\newtheorem{problem}[theorem]{Problem}
\newtheorem{algorithm}[theorem]{Algorithm}

\begin{document}
\title{Factor Analysis and Alternating Minimization}
\author{Lorenzo Finesso\footnote{Institute of Biomedical Engineering, CNR-ISIB, Padova,
\tt{lorenzo.finesso@isib.cnr.it}} \and Peter
Spreij\footnote{Korteweg-de Vries Institute for Mathematics,
Universiteit van Amsterdam, Amsterdam,
\tt{spreij@science.uva.nl}}}

\date{} \maketitle

\bigskip
\begin{center}\textit{Dedicated to Giorgio Picci on the occasion of his
65th birthday. \\ Happy Birthday Giorgio!}\end{center}

\begin{abstract}\noindent
In this paper we make a first attempt at understanding how to
build an optimal \emph{approximate} normal factor analysis model.
The criterion we have chosen to evaluate the distance between
different models is the I-divergence between the corresponding
normal laws. The algorithm that we propose for the construction of
the best approximation is of an the alternating minimization kind.
\end{abstract}

\section{Introduction}
\setcounter{equation}{0}
Factor analysis, in its original
formulation, is the linear statistical model
\begin{equation} \label{eq:fam}
Y = HX+\eps
\end{equation}
where $H$ is a deterministic matrix, $X$ and $\eps$ independent
random vectors, the first with dimension smaller than $Y$, the
second with independent components. What makes this model attractive
in applied research is the \textit{data reduction} mechanism built
in it. A large number of observed variables $Y$ are explained in
terms of a small number of unobserved (latent) variables $X$
perturbed by the independent noise $\eps$. Under normality
assumptions, which are the rule in the standard theory, all the laws
of the model are specified by covariance matrices. More precisely,
assume that $X$ and $\eps$ are zero mean independent normal vectors
with $\cov(X)=P$ and $\cov(\eps)=D$, where $D$ is diagonal. It
follows from~(\ref{eq:fam}) that $\cov(Y)= HPH^\top +D$.

Building a factor analysis model of the observed data requires the
solution of a difficult algebraic problem. Given $\Sigma_0$, the
covariance matrix of $Y$, find the triples $(H, P, D)$ such that
$\Sigma_0 = HPH^\top +D$. Due to the structural constraint on $D$,
which is assumed to be diagonal, the existence and unicity of a
factor analysis model are not guaranteed. As it turns out, the right
tools to deal with this situation come from the theory of stochastic
realization, see~\cite{finessopicci} (trying to spot the master's
hand) for an early contribution on the subject.

In the present paper we make a first attempt at understanding how to
build an optimal \emph{approximate} factor analysis model. The
criterion we have chosen to evaluate the distance between
covariances is the I-divergence between the corresponding normal
laws. The algorithm that we propose for the construction of the best
approximation is inspired by the alternating minimization procedure
of~\cite{csiszar} and~\cite{finessospreij}.

\section{The model}\setcounter{equation}{0}

Consider two independent, zero mean, normal vectors $X$ and $\eps$
of respective dimensions $k$ and $n$. We will assume that
$\cov(X)=I$, the identity matrix, and $\cov(\eps)=D>0$, a diagonal
matrix. Let $H$ be an $n \times k$ matrix (in this paper $k < n$)
and let the random vector $Y$ be defined by
\begin{equation}\label{eq:y}
Y=HX+\eps.
\end{equation}
Under these assumptions (\ref{eq:y}) is called a factor analysis
(FA) model of size $k$ for the vector $Y$. Notice that allowing
$\cov(X)=P>0$ does not produce a more general model, as a square
root of $P$ can always be absorbed in $H$. We will say that a normal
vector $Y$ admits a FA model of size $k$ if it is equal in
distribution to $HX+\eps$ for some $X$ and $\eps$ as above, i.e. if
its covariance $\Sigma_0$ can be written as $\Sigma_0 = HH^\top +D$.
Not every normal vector $Y$ admits a FA model, the hard constraint
being imposed by the diagonal structure of $D$. A probabilistic
interpretation stems from $\cov(Y|X)=D$ (see
equation~(\ref{eq:condcov}) of the Appendix) i.e. the $n$ components
of $Y$ are conditionally independent given the $k<n$ components of
some vector $X$. In Remark~\ref{rem:real} of the next section the
condition for the existence of a FA model is slightly reformulated.

\noindent Although the construction of an exact FA model is not
always possible, one can search for a best approximate model,
according to some criterion. In this paper we opt for minimizing the
I-divergence (Kullback-Leibler distance) between normal laws. Recall
that given two probability measures $\pp_1$ and $\pp_2$, defined on
the same measurable space, such that $\pp_1\ll \pp_2$, the
I-divergence of $\pp_1$ with respect to $\pp_2$ is defined as
\[
D(\pp_1||\pp_2)=\ee_{\pp_1}\log \frac{d\pp_1}{d\pp_2}.
\]
If $\pp_1$ and $\pp_2$ are normal measures on the same space $\R^n$,
with zero means and strictly positive covariance matrices $\Sigma_1$
and $\Sigma_2$ respectively, the I-divergence $ D(\pp_1||\pp_2) $
takes the explicit form
\begin{equation}\label{eq:divsigma}
D(\pp_1||\pp_2) = \half \, \log \frac{|\Sigma_2|}{|\Sigma_1|} +
\half \, {\rm tr}(\Sigma_2^{-1}\Sigma_1) - \frac{n}{2} .
\end{equation}
Since the I-divergence only depends on the covariance matrices, we
usually write $D(\Sigma_1||\Sigma_2)$ instead of $D(\pp_1||\pp_2)$.

\noindent The approximate factor analysis problem can be posed as
follows:
\begin{problem}\label{problem2}
Given the positive covariance matrix $\Sigma_0\in\openR^{n\times n}$
and the integer $k < n$  minimize
\[
D(\Sigma_0 || HH^\top+D) = \half \,\log
\frac{|HH^\top+D|}{|\Sigma_0|} +\half \,{\rm
tr}((HH^\top+D)^{-1}\Sigma_0) -\frac{n}{2}.
\]
over all pairs $(H, D)$ where $H\in\openR^{n\times k}$ and $D>0$ is
of size $n$ and diagonal.
\end{problem}

\noindent Notice that $D(\Sigma_1||\Sigma_2)$, computed as
in~(\ref{eq:divsigma}), can be considered as a divergence between
two positive definite matrices, without referring to normal
distributions. Hence Problem~\ref{problem2} also has a meaning, when
one refrains from assumptions like normality.

\noindent Existence of the minimum is guaranteed by the following
\begin{proposition}\label{prop:exist}
There exist matrices $H^*\in\openR^{n\times k}$ and $D^* >0$ of size
$n$ and diagonal minimizing the I-divergence in
Problem~\ref{problem2}.
\end{proposition}
The proof is deferred to section~\ref{section:proof}, since it uses
later results.

\noindent In order to construct an algorithm for the solution of
Problem \ref{problem2} we will imitate the approach
of~\cite{finessospreij}. The algorithm will therefore be derived by
a relaxation technique, lifting the original minimization problem to
a higher dimensional space. In the larger space a double
minimization problem equivalent to Problem~\ref{problem2} can be
formulated, leading in a natural way to an alternating minimization
algorithm.

\section{Lifting of the original problem}\setcounter{equation}{0}

In this section we will embed Problem~\ref{problem2} into a higher
dimensional space. First we introduce the relevant sets of
covariances. Given $k<n$ we denote by
\begin{equation} \label{eq:dec}
\bs=\{\Sigma \in \R^{(n+k)\times(n+k)} : \Sigma=\begin{pmatrix}
\Sigma_{11} & \Sigma_{12} \\ \Sigma_{21} & \Sigma_{22}
\end{pmatrix} >0 \}.
\end{equation}
where $\Sigma_{11}$ is $n\times n$. Two subsets of $\bs$ will play a
special role.
\begin{equation}\label{eq:dec0}
\bs_0=\{\Sigma\in\bs: \Sigma_{11}=\Sigma_0\}.
\end{equation}
where $\Sigma_0$ is a given covariance. We also consider the subset
\begin{equation} \label{eq:dec1}
\bs_1= \{\Sigma \in \bs : \Sigma=\begin{pmatrix} HH^\top +D & \quad HQ \\
(HQ)^\top & \quad Q^\top Q
\end{pmatrix} \},
\end{equation}
where $H \in \R^{n \times k}, Q \in \R^{k \times k}$ invertible,
$D>0$ diagonal. Elements of $\bs_1$ will often be denoted by
$\Sigma(H,D,Q)$.
\begin{remark} \label{rem:real} Notice that a normal vector $Y$, with $\cov(Y)=\Sigma_0$, admits a FA
model of size $k$ iff $\bs_0 \cap \bs_1 \ne \emptyset$. Supposing
that this is the case, take $\Sigma \in \bs_0 \cap \bs_1$ then, for
some $(H, D, Q)$, one has
$$\Sigma =\begin{pmatrix} \Sigma_0 & \quad HQ \\
(HQ)^\top & \quad Q^\top Q
\end{pmatrix} = \begin{pmatrix} HH^\top +D & \quad HQ \\
(HQ)^\top & \quad Q^\top Q
\end{pmatrix}>0.$$
is a \textit{bonafide} covariance of a normal vector $V$ of
dimension $n+k$. Partition $V^\top=(Y^\top, Z^\top)^\top$. It is
easy to verify that $\cov(Y) = \Sigma_0= HH^\top +D$ is the same as
$\cov(HX+\eps)$ for some $X$ standard normal and $\eps$ normal,
independent from $X$, and with diagonal covariance $D$.
\end{remark}

\medskip
\noindent The lifted  minimization problem can be posed as follows
\begin{problem}\label{liftedproblem}
\[
\min_{\Sigma'\in\bs_0,\Sigma_1\in\bs_1}D(\Sigma'||\Sigma_1)
\]
\end{problem}

\noindent which can be viewed as an iterated minimization problem
over each of the variables. The two resulting partial minimization
problems will be investigated in the following sections. In
section~\ref{link} we will show the connection between
Problems~\ref{problem2} and~\ref{liftedproblem}. More precisely, we
will prove
\begin{proposition}\label{prop:pqq}
Let $\Sigma_0$ be given. It holds that
\begin{equation*}
\min_{H,D}D(\Sigma_0||HH^\top +
D)=\min_{\Sigma'\in\bs_0,\Sigma_1\in\bs_1}D(\Sigma'|\Sigma_1).
\end{equation*}
\end{proposition}

\subsection{The first partial minimization problem}

In this section we consider the first of the two partial
minimization problems. Here we minimize, for a given positive
definite matrix $\Sigma\in \bs$, the divergence $D(\Sigma'||\Sigma)$
over $\Sigma'\in\bs_0$. The unique solution to this problem can be
computed analytically and follows from
\begin{lemma}\label{lemma:pm1}
Let $(Y,X)$ be a random vector distributed according to some
$Q=Q^{Y,X}$ and let $\bp$ the set of all distributions $P=P^{Y,X}$
whose marginal $P^Y=P_0$, for some fixed $P_0 \ll Q^Y$. Then
$\min_{P\in \bp} D(P||Q)= D(P^*||Q)$ where $P^*$ is given by the
Radon-Nikodym derivative
\[
\frac{dP^*}{dQ}=\frac{dP_0}{dQ^Y}.
\]
Moreover,
\begin{equation}\label{eq:pq0}
D(P^*||Q)=D(P_0||Q^Y).
\end{equation}
and, for any other $P \in \bp$, one has the Pythagorean law
\begin{equation}\label{eq:p1}
D(P||Q)=D(P||P^*)+D(P^*||Q).
\end{equation}
\end{lemma}
\begin{proof}
First we show that~(\ref{eq:pq0}) holds. Recall that $Y$ has law
$P_0$ under $P^*$, then
\[
D(P^*||Q)= \ee_{P^*}\log \frac{dP^*}{dQ} =  \ee_{P^*}\log
\frac{dP_0}{dQ^Y} = \ee_{P_0}\log \frac{dP_0}{dQ^Y} = D(P_0||Q^Y).
\]
To show that $P^*$ is a minimizer it is clearly sufficient to prove
that~(\ref{eq:p1}) holds.
\begin{align*}
D(P||Q) & =\ee_P\log \frac{dP}{dP^*}+\ee_P\log \frac{dP^*}{dQ} \\
& = D(P||P^*) + \ee_P\log \frac{dP_0}{dQ^Y} \\
& = D(P||P^*) + \ee_{P_0}\log \frac{dP_0}{dQ^Y} =
D(P||P^*)+D(P^*||Q),
\end{align*}
where we used the fact that all $P \in \bp$ have marginal $P^Y=P_0$.
\end{proof}
\begin{remark}\label{remark:pm1}
The law $P^*$ is easily characterized in terms of the problem data
$P_0$ and $Q$ noticing that the marginal $P^{* Y}=P_0$ and the
conditional $P^{* X|Y}=Q^{X|Y}$.
\end{remark}
We now apply Lemma~\ref{lemma:pm1} to the case of normal laws and
solve the first partial minimization. See also~\cite{cramer1} for a
different proof.
\begin{proposition}\label{prop_min_first}
Let $Q$ and $P_0$ be zero mean normal laws with strictly
po\-si\-ti\-ve covariances $\Sigma \in \bs$ and $\Sigma_0 \in \R^{n
\times n}$ respectively. Then, $\min_{\Sigma' \in \bs_0} D(\Sigma'
|| \Sigma)$ is attained by the zero mean normal law $P^*$ with
covariance
\begin{equation}\nonumber \Sigma^*=\begin{pmatrix} \Sigma_0 & \quad
\Sigma_0\Sigma_{11}^{-1}\Sigma_{12}
\\
\Sigma_{21}\Sigma_{11}^{-1}\Sigma_0 & \quad
\Sigma_{22}-\Sigma_{21}\Sigma_{11}^{-1}
(\Sigma_{11}-\Sigma_{0})\Sigma_{11}^{-1}\Sigma_{12}
\end{pmatrix} >0.
\end{equation}
Moreover,
$$D(\Sigma^*||\Sigma)= D(\Sigma_0||\Sigma_{11}).$$
\end{proposition}
\begin{proof} This follows from
Remark~\ref{remark:pm1}. A direct computation gives
\begin{align*}\Sigma^*_{12} = & \,\, \mathbb{E}_{P^*}XY^\top = \,\,
\mathbb{E}_{P^*}(\mathbb{E}_{P^*}[X|Y]Y^\top)
\\
 = & \,\, \mathbb{E}_{P^*}(\mathbb{E}_Q[X|Y]Y^\top)
 = \,\, \mathbb{E}_{P^*}(\Sigma_{21}\Sigma_{11}^{-1}YY^\top) \\
 = & \,\, \Sigma_{21}\Sigma_{11}^{-1}\,\mathbb{E}_{P_0}YY^\top
 \,\, = \,\, \Sigma_{21}\Sigma_{11}^{-1}\Sigma_0.
\end{align*}
Likewise, we have
\begin{align*}\Sigma^*_{22} = & \,\, \mathbb{E}_{P^*}XX^\top  = \cov_{P^*}(X) \\
& = \cov_{P^*}(X|Y) + \mathbb{E}_{P^*}
(\mathbb{E}_{P^*}[X|Y]\mathbb{E}_{P^*}[X|Y]^\top) \\
& = \cov_{Q}(X|Y) + \mathbb{E}_{P^*}
(\mathbb{E}_{Q}[X|Y]\mathbb{E}_{Q}[X|Y]^\top) \\
& = \Sigma_{22}-\Sigma_{21}\Sigma_{11}^{-1}\Sigma_{12} +
\mathbb{E}_{P^*}
(\Sigma_{21}\Sigma_{11}^{-1}Y(\Sigma_{21}\Sigma_{11}^{-1}Y)^\top) \\
& = \Sigma_{22}-\Sigma_{21}\Sigma_{11}^{-1}\Sigma_{12} +
\mathbb{E}_{P_0}
(\Sigma_{21}\Sigma_{11}^{-1}YY^\top\Sigma_{11}^{-1}\Sigma_{12}) \\
& = \Sigma_{22}-\Sigma_{21}\Sigma_{11}^{-1}\Sigma_{12} +
\Sigma_{21}\Sigma_{11}^{-1}\Sigma_0\Sigma_{11}^{-1}\Sigma_{12}.
\end{align*}
Notice that, since $\Sigma>0$ by assumption,
\[
\Sigma^*_{22}-\Sigma^*_{21}
(\Sigma^*_{11})^{-1}\Sigma^*_{12}=\Sigma_{22}-\Sigma_{21}\Sigma_{11}^{-1}\Sigma_{12}
>0
\]
which, together with the assumption $\Sigma_0>0$, shows that
$\Sigma^*>0$. \\ \noindent The last relation,
$D(\Sigma^*||\Sigma)=D(\Sigma_0||\Sigma_{11})$, reflects
equation~(\ref{eq:pq0}).
\end{proof}

\subsection{The second partial minimization problem}

In this section we turn to the second partial minimization problem.
Here we minimize, for a given positive definite matrix $\Sigma\in
\bs$, the divergence $D(\Sigma||\Sigma_1)$ over $\Sigma_1\in\bs_1$.

\noindent Clearly this problem cannot have a unique solution in
terms of the matrices $H$ and $Q$. Indeed, if $U$ is a unitary
$k\times k$ matrix and $H'=HU$, $Q'=U^\top Q$, then
$H'H'^\top=HH^\top$, $Q'^\top Q'=Q^\top Q$ and $H'Q'=HQ$.
Nevertheless, the optimal matrices $HH^\top$, $HQ$ and $Q^\top Q$
are unique, as we will see in Proposition~\ref{prop:pm2}. First we
need to introduce some notation and conventions. If $P$ is a
positive definite matrix, we denote by $P^{1/2}$ any matrix
satisfying $(P^{1/2})^\top (P^{1/2})=P$, and by $P^{-1/2}$ its
inverse. If $M$ is any square matrix, we denote by $\Delta(M)$ the
diagonal matrix
\[
\Delta(M)_{ii}=M_{ii}.
\]
Recall that  we denote by $\Sigma(H,D,Q)$ a typical element of
$\bs_1$.
\begin{proposition}\label{prop:pm2}
Given $\Sigma\in\bs$ the $\min_{\Sigma_1 \in \bs_1 } D(\Sigma ||
\Sigma_1)$ is attained at a $\Sigma_1^*$ such that
$\Sigma_1\in\bs_1$ is solved by
\begin{align*}
Q^* & =\Sigma_{22}^{1/2},\\
H^* & = \Sigma_{12}\Sigma_{22}^{-1/2}, \\
D^* & = \Delta(\Sigma_{11}-\Sigma_{12}\Sigma_{22}^{-1}\Sigma_{21}).
\end{align*}
Thus the minimizing matrix $\Sigma_1^*=\Sigma(H^*,D^*,Q^*)$ becomes
\begin{equation} \label{eq:sigmastar}
\Sigma_1^*=\begin{pmatrix}
\Sigma_{12}\Sigma_{22}^{-1}\Sigma_{21}+\Delta(\Sigma_{11}-\Sigma_{12}\Sigma_{22}^{-1}\Sigma_{21})
& \quad \Sigma_{12} \\
\Sigma_{21} & \quad \Sigma_{22}
\end{pmatrix}.
\end{equation}
Moreover, the Pythagorean law
\begin{equation}\label{eq:p2}
D(\Sigma||\Sigma(H,D,Q))=D(\Sigma||\Sigma_1^*)+D(\Sigma_1^*||\Sigma(H,D,Q))
\end{equation}
holds for any $\Sigma(H,D,Q)\in\bs_1$, and therefore  $\Sigma_1^*$
is unique.
\end{proposition}
\begin{proof}
It is sufficient to show the validity of~(\ref{eq:p2}). We first
compute
$$2 D(\Sigma||\Sigma(H,D,Q))- 2 D(\Sigma||\Sigma^*_1).$$
It follows from Lemma~A.\ref{lemma:matrix} that
$|\Sigma(H,D,Q)|=|D|\times |Q^\top Q|$. In view of
equation~(\ref{eq:divsigma}) the above difference becomes
\begin{equation} \label{eq:tr3}
\log |D| + \log |Q^\top Q| - \log |D^*| - \log |Q^{*^\top} Q^*|+
{\rm tr} \big( \Sigma(H,D,Q)^{-1}\Sigma\big)-{\rm tr} \big(
\Sigma^{*-1}_1\Sigma\big).
\end{equation}
Using Corollary~A.\ref{cor:inv}, we compute
\begin{equation}\label{eq:sinvs}
\Sigma(H,D,Q)^{-1}=
\begin{pmatrix}
D^{-1} & \quad  -D^{-1}HQ^{-\top} \\
-Q^{-1}H^\top D^{-1} & \quad Q^{-1}(H^\top D^{-1}H+ I)Q^{-\top}
\end{pmatrix},
\end{equation}
and hence we get that
\begin{align}\label{eq:tr1}
\lefteqn{{\rm tr}\big(\Sigma(H,D,Q)^{-1}\Sigma\big)    = {\rm
tr}\big(D^{-1}(\Sigma_{11}-HQ^{-\top}\Sigma_{21}) \big)} \nonumber\\
& ~\mbox{} + {\rm tr}\big(-Q^{-1}H^\top
D^{-1}\Sigma_{12}+Q^{-1}(H^\top D^{-1}H+
I)Q^{-\top}\Sigma_{22} \big) \nonumber\\
& = {\rm
tr}\big(D^{-1}(\Sigma_{11}-2HQ^{-\top}\Sigma_{21})+Q^{-1}(H^\top
D^{-1}H+ I)Q^{-\top}\Sigma_{22} \big).
\end{align}
Apply now Lemma~A.\ref{lemma:matrix} to~(\ref{eq:sigmastar}) and
write $\Delta=
\Delta(\Sigma_{11}-\Sigma_{12}\Sigma_{22}^{-1}\Sigma_{21})$, to get
\[
\Sigma^{*-1}_1 = \Sigma(H^*,D^*,Q^*)^{-1}=
\begin{pmatrix}
\Delta^{-1} & \quad -\Delta^{-1}\Sigma_{12}\Sigma_{22}^{-1} \\
-\Sigma_{22}^{-1}\Sigma_{21}\Delta^{-1} & \quad
\Sigma_{22}^{-1}\Sigma_{21}\Delta^{-1}\Sigma_{12}\Sigma_{22}^{-1}+\Sigma_{22}^{-1}
\end{pmatrix}.
\]
Therefore
\begin{equation}\label{eq:tr2}
{\rm tr}\big(\Sigma^{*-1}_1\Sigma\big) = {\rm tr}\big(\Delta^{-1}
\times
(\Sigma_{11}-\Sigma_{12}\Sigma_{22}^{-1}\Sigma_{21})\big)+{\rm
tr}\,I_k = {\rm tr}\big(\Delta^{-1}\Delta) + k = n+k.
\end{equation}
Combining equations~(\ref{eq:tr3}), (\ref{eq:tr1}),
and~(\ref{eq:tr2}), we find that
\begin{align}\label{eq:divdiff}
\lefteqn{D(\Sigma||\Sigma(H,D,Q))-D(\Sigma||\Sigma^*_1) =}
\nonumber\\
& \qquad \log |D| + \log |Q^\top Q| - \log |D^*| - \log |Q^{*^\top}
Q^*|
\nonumber\\
& \qquad + \mbox{} {\rm
tr}\big(D^{-1}(\Sigma_{11}-HQ^{-\top}\Sigma_{21})
\big) \nonumber\\
 & \qquad + {\rm tr}\big(-Q^{-1}H^\top
D^{-1}\Sigma_{12}+Q^{-1}(H^\top D^{-1}H+
I)Q^{-\top}\Sigma_{22}\big) \nonumber \\
& \qquad -(n+k).
\end{align}
We proceed with the computation of $2 D(\Sigma^*_1||\Sigma(H,D,Q))$.
\begin{align}\label{eq:starss}
\lefteqn{2 D(\Sigma(H^*,D^*,Q^*)||\Sigma(H,D,Q))=}\nonumber\\
& \qquad \log |D|+\log |Q^\top Q| - \log |D^*|-\log |Q^{*^\top}
Q^*|-(n+k) \nonumber\\
& \qquad +{\rm tr}\big(\Sigma(H,D,Q)^{-1}\Sigma(H^*,D^*,Q^*)\big).
\end{align}
Combining equations~(\ref{eq:sigmastar}), (\ref{eq:sinvs}),  and
${\rm tr}\big(
D^{-1}(\Sigma_{12}\Sigma_{22}^{-1}\Sigma_{21}+\Delta)\big)={\rm
tr}\big( D^{-1}\Sigma_{11}\big)$, we obtain
\begin{align}\label{eq:laststep}
& {\rm tr}\big(\Sigma(H,D,Q)^{-1}\Sigma(H^*,D^*,Q^*)\big)= {\rm tr}\big( D^{-1}\Sigma_{11}\big) \nonumber \\
& \quad -2 {\rm tr}\big(D^{-1}HQ^{-\top}\Sigma_{21}\big)  + {\rm
tr}\big(Q^{-1}(H^\top D^{-1}H+ I)Q^{-\top}\Sigma_{22} \big).
\end{align}
Insertion of~(\ref{eq:laststep}) into~(\ref{eq:starss}) and a
comparison with~(\ref{eq:divdiff}) yields the result.
\end{proof}
\begin{remark}\label{remark:hhd}
Notice that the matrix $H^*H^{*\!^\top}$ is strictly dominated by
$\Sigma_{11}$ (in the sense of positive matrices). This easily
follows from
$\Sigma_{11}-H^*H^{*^\top}=\Sigma_{11}-\Sigma_{12}\Sigma_{22}^{-1}\Sigma_{21}>0$,
and the assumption $\Sigma>0$. By the same token $D^*>0$.
\end{remark}

\subsection{The link to the original
problem}\label{link}

We now establish the connection between the lifted problem and the
original Problem~\ref{problem2}. \medskip\\
{\bf Proof of Proposition~\ref{prop:pqq}} Let
$\Sigma_1=\Sigma(H,D,Q)$ and denote by
$\Sigma^*=\Sigma^*(\Sigma_1)$, the solution of the first partial
minimization over $\bs_0$. We have, for all $\Sigma' \in \bs_0$,
\begin{align*}
D(\Sigma'||\Sigma_1)& \geq D(\Sigma^*||\Sigma_1) \\
& = D(\Sigma_0||HH^\top + D) \\
& \geq  \min_{H,D}D(\Sigma_0||HH^\top + D),
\end{align*}
where we used Proposition~\ref{prop:exist} to write $\min$ on the
RHS. It follows that
\[
\inf_{\Sigma'\in\bs_0,\Sigma_1\in\bs_1}D(\Sigma' ||
\Sigma_1)\geq\min_{H,D}D(\Sigma_0||HH^\top + D).
\]
Conversely, let $(H^*,D^*)$ be the minimizer of $(H,D)\mapsto
D(\Sigma_0||HH^\top+D)$, pick an arbitrary invertible $Q^*$, and let
$\Sigma^*=\Sigma(H^*,D^*,Q^*)$ be the corresponding element in
$\bs_1$. Furthermore, let $\Sigma^{**}\in\bs_0$ be the minimizer of
$\Sigma\mapsto D(\Sigma||\Sigma^*)$ over $\bs_0$. Then
\begin{align*}
\min_{H,D}D(\Sigma_0||HH^\top + D) = & \,\, D(\Sigma_0||H^*H^{*^\top}+D^*) \\
\geq  & \,\, D(\Sigma^{**}|| \Sigma^*) \\
\geq & \inf_{\Sigma'\in\bs_0,\Sigma_1\in\bs_1}D(\Sigma'||\Sigma_1),
\end{align*}
which shows the opposite inequality. Finally, to show that we can
replace the infima with minima also in the lifted problem, notice
that (see Proposition~\ref{prop_min_first}) $D(\Sigma^{**}||
\Sigma^*)=D(\Sigma_0||H^*H^{*^\top}+D^*)$.

\section{Alternating minimization algorithm}\setcounter{equation}{0}

In this section we combine the two partial minimization problems
above to derive an iterative algorithm for Problem~\ref{problem2}.
It turns out that this algorithm is also instrumental in proving the
existence of a solution to Problem~\ref{problem2}.

\subsection{The algorithm}

We suppose that the given matrix $\Sigma_0$ is strictly positive
definite. Pick the initial values $H_0, D_0, Q_0$ such that $H_0$ is
of full rank, $D_0>0$ is diagonal, $Q_0$ and $H_0H_0^\top +D_0$ are
invertible.
\medskip\\
At the $t$-th iteration the matrices $H_t$, $D_t$ and $Q_t$ are
available. Start solving the first partial minimization problem with
$\Sigma = \Sigma(H_t,D_t,Q_t)$. Use the resulting matrix as data for
the second partial minimization, the solution of which gives the
update rules
\begin{align}
Q_{t+1} & = \Big(Q_t^\top Q_t - Q_t^\top H_t^\top (H_tH_t^\top +
D_t)^{-1}H_tQ_t \nonumber \\ & \quad + Q_t^\top H_t^\top
(H_tH_t^\top + D_t)^{-1}\Sigma_0 (H_tH_t^\top +
D_t)^{-1}H_tQ_t\Big)^{1/2}, \label{eq:qn} \\
H_{t+1} & = \Sigma_0(H_tH_t^\top + D_t)^{-1}H_tQ_tQ_{t+1}^{-1}, \label{eq:h1}\\
D_{t+1} & = \Delta(\Sigma_0-H_{t+1}H_{t+1}^\top).\label{eq:dn}
\end{align}
In~(\ref{eq:qn}) there is some freedom in computing the square root
that determines $Q_{t+1}$. Properly choosing the square root will
result in the disappearance of $Q_t$ from the algorithm. This is an
attractive feature, since $Q_t$ only serves as an auxiliary
variable. One can write the RHS of equation~(\ref{eq:qn}), before
taking the square root, as
\[
Q_t^\top(I- H_t^\top (H_tH_t^\top + D_t)^{-1}(H_tH_t^\top + D_t -
\Sigma_0) (H_tH_t^\top + D_t)^{-1}H_t)Q_t
\]
and denoting
\begin{equation} \label{eq:rt}
R_t = I- H_t^\top (H_tH_t^\top + D_t)^{-1}(H_tH_t^\top + D_t -
\Sigma_0) (H_tH_t^\top + D_t)^{-1}H_t
\end{equation}
a possible square root is given by
\[
R_t^{1/2}Q_t.
\]
Notice that $R_t$ only involves the iterates $H_t$ and $D_t$. The
update equation~(\ref{eq:qn}) can therefore be rewritten as
\begin{equation}\label{eq:hn}
H_{t+1}  = \Sigma_0(H_tH_t^\top + D_t)^{-1}H_tR_t^{-1/2}.
\end{equation}
The final version of the algorithm is given by
equations~(\ref{eq:dn}),(\ref{eq:rt}), and~(\ref{eq:hn})  which, for
clarity, we present as
\begin{algorithm}\label{algo1}
\begin{align}
R_t & = I- H_t^\top (H_tH_t^\top + D_t)^{-1}(H_tH_t^\top + D_t -
\Sigma_0) (H_tH_t^\top + D_t)^{-1}H_t, \label{eq:r} \\
H_{t+1} & = \Sigma_0(H_tH_t^\top + D_t)^{-1}H_tR_t^{-1/2}, \label{eq:h}\\
D_{t+1} & = \Delta(\Sigma_0-H_{t+1}H_{t+1}^\top). \label{eq:d}
\end{align}
\end{algorithm}
\noindent In order to avoid taking a square root at each step one
can introduce the matrices $K_t=H_tQ_t$ and $P_t = Q_t^TQ_t$ and
write the updates for $K_t$ and $P_t$. Equations~(\ref{eq:qn}),
(\ref{eq:h1}), and (\ref{eq:dn}) easily give
\begin{algorithm}\label{algo2}
\begin{align}
K_{t+1} & =\Sigma_0(K_tP_t^{-1}K_t^\top + D_t)^{-1}K_t, \label{eq:k}\\
P_{t+1} & = P_t-K_t^\top (K_tP_t^{-1}K_t^\top+ D_t)^{-1}
(K_tP_t^{-1}K_t^\top-\Sigma_0)(K_tP_t^{-1}K_t^\top+D_t)^{-1}K_t, \nonumber\\
D_{t+1} & =
\Delta(\Sigma_0-K_{t+1}P_{t+1}^{-1}K_{t+1}^\top).\nonumber
\end{align}
After the final iteration, the $T$-th say, one can take
$H_T=K_TQ_T^{-1}$, where $Q_T$ is a square root of $P_T$.
\end{algorithm}

\noindent Notice that in both Algorithm~\ref{algo1} and~\ref{algo2}
it is required to invert $n\times n$ matrices (like e.g.
$(H_tH_t^\top + D_t)^{-1}$). Applying corollary~A.\ref{cor:inv} one
gets $(H_tH_t^\top + D_t)^{-1}H_t=D_t^{-1}H_t(I+H_t^\top
D_t^{-1}H_t)$. Hence, we can replace e.g.~(\ref{eq:hn}) with
\begin{equation}\label{eq:ht}
H_{t+1}  = \Sigma_0D_t^{-1}H_t(I+H_t^\top
D_t^{-1}H_t)^{-1}R_t^{-1/2}.
\end{equation}
By the same token one can write $$K_{t+1}  =
\Sigma_0D_t^{-1}K_t(P_t+K_t^\top D_t^{-1}K_t)^{-1}P_t$$ to replace
~(\ref{eq:k}).

\noindent Some properties of the algorithm are summarized in the
next proposition.
\begin{proposition}\label{prop:properties}
For Algorithm~\ref{algo1} the following hold for all $t$. \\
(a) $D_t>0$ and $(D_t)_{ii} \le (\Sigma_0)_{ii}$.
\\
(b) $R_t$ is invertible. \\
(c) If $H_0$ is of full column rank, so is $H_t$. \\
(d) $H_tH_t^\top \leq \Sigma_0$. \\
(e) If $\Sigma_0 = H_0H_0^\top+D_0$ then the algorithm stops. \\
(f) The objective function decreases at each iteration. More
precisely, let $\Sigma_{0,t}$ be the solution of the first partial
minimization with data $\Sigma_t=\Sigma(H_t,D_t,Q_t)$. Then
\begin{equation}\nonumber
D(\Sigma_0||H_{t+1}H_{t+1}^\top) - D(\Sigma_0|| H_{t}H_{t}^\top) =
-\Big(D(\Sigma_{t+1}||\Sigma_t)+D(\Sigma_{0,t}||
\Sigma_{0,t+1})\Big).
\end{equation}
(g) The limit points $(H,D)$ of the algorithm satisfy the relations
\begin{align*}
H  & = (\Sigma_0-HH^\top)D^{-1}H, \\
D & = \Delta(\Sigma_0-HH^\top).
\end{align*}
\end{proposition}

\begin{proof}
(a) This follows from Remark~\ref{remark:hhd}.\\
(b) Use the identity $I-H_t^\top(H_tH_t^\top +
D_t)^{-1}H_t=(I+H_t^\top D_t^{-1}H_t)^{-1}$ and the assumption
$\Sigma_0>0$. \\
(c) Use the assumption $\Sigma_0>0$, (a), and (b).\\
(d) Again from Remark~\ref{remark:hhd} and the construction of the
algorithm as a combination of the two partial minimization
problems.  \\
(e) This is a triviality upon noticing that one can take $R_t=I$ in
this
case. \\
(f) It follows from a concatenation of Lemma~\ref{lemma:pm1} and
Proposition~\ref{prop:pm2}. Notice that we can express the decrease
as the sum of two I-divergences, since the Pythagorean law holds for
both partial minimizations.  \\
(g) We consider Algorithm~\ref{algo2} first. Assume that all
variables converge. Then, from~(\ref{eq:k}), the limit points $K,
P, D$ satisfy the relation $K=\Sigma_0D^{-1}K(P+K^\top
D^{-1}K)^{-1}P$. Postmultiplication  by $P^{-1}(P+K^\top D^{-1}K)$
yields, after rearranging terms,
$K=(\Sigma_0-KP^{-1}K^\top)D^{-1}K$. Let now $Q$ be a square root
of $P$ and $H=KQ^{-1}$ to get the first relation. The rest is
trivial.
\end{proof}

\subsection{Proof of
Proposition~\ref{prop:exist}}\label{section:proof}

Let $D_0$ and $H_0$ be arbitrary and perform one step of the
algorithm to get matrices $D_1$ and $H_1$. It follows from
Proposition~\ref{prop:properties} that $D(\Sigma_0||H_1H_1^\top
+D_1)\leq D(\Sigma_0||H_0H_0^\top +D_0)$. Moreover, $H_1H_1^\top
\leq \Sigma_0$ and $D_1\leq \Delta(\Sigma_0)$. Hence the search
for a minimum can be confined to the set of matrices $(H, D)$
satisfying $HH^\top \leq \Sigma_0$ and $D\leq \Delta(\Sigma_0)$.
Next, we claim that it is also sufficient to restrict the search
for a minimum to all matrices $(H, D)$ such that $HH^\top +D\geq
\eps I$ for some sufficiently small $\eps>0$. Indeed, if the last
inequality is violated, then $HH^\top +D$ has an eigenvalue less
than $\eps$. Write the Jordan decompositions $HH^\top + D=U\Lambda
U^\top$, and let $\Sigma_U=U^\top\Sigma_0U$. Then
$D(\Sigma_0||HH^\top + D)=D(\Sigma_U||\Lambda)$, as one easily
verifies. Denoting by $\lambda_i$ the eigenvalues of $HH^\top + D$
and letting $\sigma_{ii}$ be the diagonal elements of $\Sigma_U$,
we can write $D(\Sigma_U| \Lambda)=-\half \log
|\Sigma_U|+\half\sum_i\log\lambda_i-
\frac{n}{2}+\half\sum_i\frac{\sigma_{ii}}{\lambda_i}$. Let
$\lambda_{i_0}$ be a minimum eigenvalue and take $\eps$ smaller
than the minimum of all $\sigma_{ii}$, which is positive, since
$\Sigma_0$ is strictly positive definite. Then the contribution
for $i=i_0$ in the summation to the divergence
$D(\Sigma_U||\Lambda)$ is at least $\log\eps +1$, which tends to
infinity for $\eps\to 0$. This proves the claim. So, we have shown
that a minimizing pair $(H,D)$ has to satisfy $HH^\top\leq
\Sigma_0$, $D\leq \Delta(\Sigma_0)$, and $HH^\top + D\geq \eps I$,
for some $\eps > 0$. In other words we have to minimize the
I-divergence over a compact set on which it is clearly continuous.
This proves Proposition~\ref{prop:exist}.

\bigskip \bigskip
\appendix
\section{Appendix}\setcounter{equation}{0}

For ease of reference we collect here some standard formulas for the
normal distribution and some matrix algebra.

\subsection{Multivariate normal distribution}\label{gauss}

Let $(X^\top,Y^\top)^\top$ be a zero mean normal vector with
covariance matrix
\[
\Sigma=\begin{pmatrix} \Sigma_{XX} & \Sigma_{XY} \\
\Sigma_{YX} & \Sigma_{YY}
\end{pmatrix}.
\]
Assume that $\Sigma_{YY}$ is invertible. The conditional law of $X$
given $Y$ is normal with $\ee[X|Y] =\Sigma_{XY}\Sigma_{YY}^{-1}Y$
and
\begin{equation}\label{eq:condcov}
\cov[X|Y]=\Sigma_{XX}- \Sigma_{XY}\Sigma_{YY}^{-1}\Sigma_{YX}.
\end{equation}

\subsection{Partitioned matrices}

\begin{lemma}\label{lemma:matrix}
Let $A, D$ be square matrices. Assume invertibility where required.
\[
\begin{pmatrix}
A & C \\B & D
\end{pmatrix}=
\begin{pmatrix}
I & \quad CD^{-1} \\ 0 & \quad I
\end{pmatrix}
\begin{pmatrix}
A-CD^{-1}B & \quad 0 \\ 0 & D
\end{pmatrix}
\begin{pmatrix}
I & 0 \\D^{-1}B & \quad I
\end{pmatrix},
\]
\[
\begin{pmatrix}
A & C \\B & D
\end{pmatrix}=
\begin{pmatrix}
I & \quad 0 \\ BA^{-1} & \quad I
\end{pmatrix}
\begin{pmatrix}
A & \quad 0 \\ 0 & \quad D-BA^{-1}C
\end{pmatrix}
\begin{pmatrix}
I & \quad A^{-1}C \\0 & \quad I
\end{pmatrix},
\]
\begin{align*}
& \begin{pmatrix} A & C \\B & D
\end{pmatrix}^{-1}= \\
& \quad \begin{pmatrix} (A-CD^{-1}B)^{-1} & \quad
-(A-CD^{-1}B)^{-1}CD^{-1}
\\-D^{-1}B(A-CD^{-1}B)^{-1} &
\quad D^{-1}B(A-CD^{-1}B)^{-1}CD^{-1}+D^{-1}
\end{pmatrix}.
\end{align*}
\end{lemma}

\medskip
\begin{corollary}\label{cor:inv}
\[
(D-BAC)^{-1}=D^{-1}+D^{-1}B(A^{-1}-CD^{-1}B)^{-1}CD^{-1}.
\]
\end{corollary}
\begin{proof}
For Lemma~\ref{lemma:matrix} a check will suffice. The Corollary
follows using the two decompositions of the Lemma with $A$
replaced by $A^{-1}$ and comparing the two expressions of the
lower right block of the inverse matrix.
\end{proof}

%%%%%%%%%%%%%%%%%%%%%%%%%%%%%%%%%%%%%%%%%%%%%%%%%%%%%%%%%%%%%%%%%%%%%%

%%%%%%%%%%%%%%%%%%%%%%%%%%%%%%%%%%%%%%%%%%%%%%%%%%%%%%%%%%%%%%%%%%%%%%

\end{document}